\documentclass[10pt,conference,letterpaper]{IEEEtran}
\IEEEoverridecommandlockouts

\usepackage{cite}
\usepackage{graphicx}
\usepackage[cmex10]{amsmath}
\usepackage{amssymb}
\newcommand{\set}[1]{ {\mathcal{#1}} }
\newcommand{\vect}[1]{ {\boldsymbol{#1}} }
\newcommand{\expect}[1]{ {\mathbb{E}\left[ {#1} \right]} }
\newcommand{\prob}[1]{ {\mathbb{P}\left\{ {#1} \right\}} }

\newcommand{\abs}[1]{ {\left| {#1} \right|} }

\newcommand{\norm}[1]{ {\left\lVert {#1} \right\rVert} }
\newcommand{\indicator}[1]{ {\mathbb{I}\left\{ {#1} \right\}} }

\newcommand{\RealSet}{ {\mathbb{R}} }
\newcommand{\IntSet}{ {\mathbb{Z}} }
\newcommand{\defeq}{ {\triangleq} }

\DeclareMathOperator{\argmax}{argmax}
\DeclareMathOperator{\arginf}{arginf}

\newcommand{\maximize}{ {\text{Maximize}} }
\newcommand{\minimize}{ {\text{Minimize}} }
\newcommand{\subjectto}{ {\text{Subject to}} }

\newcommand{\prts}[1]{ {\left[ {#1} \right]} }
\newcommand{\prtr}[1]{ {\left( {#1} \right)} }
\newcommand{\prtc}[1]{ {\left\{ {#1} \right\}} }

\newcommand{\maxvar}{ {(\text{max})} }

\newcommand{\optvar}{ {(\text{opt})} }

\newcommand{\dppvar}{ {(\text{DPP})} }
\newcommand{\fqvar}{ {(\text{FQ})} }
\newcommand{\dvar}{ {\gamma} }
\newcommand{\spath}{ {\Theta} }
\newcommand{\HH}{ {\text{H}} }
\newcommand{\LL}{ {\text{L}} }
\newcommand{\KK}{ {\text{A}} }

\newtheorem{theorem}{Theorem}
\newtheorem{lemma}{Lemma}

\usepackage{algorithmic}
\newsavebox{\ieeealgbox}
\newenvironment{boxedalgorithmic}
  {\begin{lrbox}{\ieeealgbox}
   \begin{minipage}{\dimexpr\columnwidth-2\fboxsep-2\fboxrule}
   \begin{algorithmic}}
  {\end{algorithmic}
   \end{minipage}
   \end{lrbox}\noindent\fbox{\usebox{\ieeealgbox}}}

\begin{document}

\title{Achieving Utility-Delay-Reliability Tradeoff in Stochastic Network Optimization with Finite Buffers}

\author{
  \IEEEauthorblockN{Sucha Supittayapornpong,~~Michael J. Neely}
  \IEEEauthorblockA{Department of Electrical Engineering\\
    University of Southern California\\
    Email: supittay@usc.edu,~~mjneely@usc.edu}

  \thanks{This work is supported by the NSF under Career grant CCF-0747525.}
}


\maketitle

\begin{abstract}
One practical open problem is the development of a distributed algorithm that achieves near-optimal utility using only a finite (and small) buffer size for queues in a stochastic network.  This paper studies utility maximization (or cost minimization) 
in a finite-buffer regime and considers the corresponding delay and reliability (or rate of packet drops) tradeoff. A \emph{floating-queue} algorithm allows the stochastic network optimization framework to be implemented with finite buffers at the cost of packet drops.  Further, the buffer size requirement is significantly smaller than previous works in this area.  With a finite buffer size of $B$ packets, the proposed algorithm achieves within $O(e^{-B})$ of the optimal utility while maintaining average per-hop delay of $O(B)$ and an average per-hop drop rate of $O(e^{-B})$ in steady state.  From an implementation perspective, the floating-queue algorithm requires little modification of the well-known Drift-Plus-Penalty policy (including MaxWeight and Backpressure policies).  As a result, the floating-queue algorithm inherits the distributed and low complexity nature of these policies.
\end{abstract}

\section{Introduction}
Stochastic network optimization is a general framework for solving a network optimization problem with randomness \cite{Neely:SNO}.  The framework generates a control algorithm that achieves a specified objective, such as minimizing power cost or maximizing throughput utility.  It is assumed that the network has random states that evolve over discrete time.  Every time slot, a network controller observes the current network state and makes a control decision. The network state and control decision together incur some cost and, at the same time, serves some amount of traffic from network queues.  The algorithm is designed to greedily minimize a \emph{drift-plus-penalty} expression every slot. This greedy procedure is known to minimize time average network cost subject to queue stability.

This general framework has been used to solve several network optimization problems such as network routing \cite{Tassiulas:Stability}, throughput maximization \cite{Eryilmaz:QLB}, dynamic power allocation \cite{Neely:Power}, quality of information maximization \cite{Sucha:QoI}.  The framework yields low-complexity algorithms which do not require any statistical knowledge of the network states.  Therefore, these algorithms are easy to implement and are robust to environment changes.  Further, they achieve an $[O(1/V), O(V)]$ utility-delay tradeoff, where $V > 0$ is a parameter that can be chosen as desired to achieve a specific operating point on the $[O(1/V), O(V)]$ utility-delay tradeoff curve.

Prior works attempt to improve network delay without sacrificing reliability, where reliability is measured by the rate of packet drops.  Previous works \cite{Neely:SuperFast} and \cite{Neely:EnergyDelay} use an exponential Lyapunov function and assumed knowledge of an $\epsilon$ parameter, where $\epsilon$ measures a distance associated with the optimal operation point.  They achieve an optimal $[O(1/V), O(\log(V))]$ utility-delay tradeoff.   A simpler methodology allows packet drops in order to obtain an $[O(1/V), O([\log(V)]^2)]$ utility-delay tradeoff \cite{Longbo:Lagrange,Longbo:LIFO}.   In \cite{Longbo:Lagrange}, a steady state behavior is observed to learn a placeholder parameter to achieve the tradeoff in steady state.  However, the algorithm does not gracefully adapt to changes of the network state distribution.  It would 
need another mechanism to sense changes and then recompute a new placeholder parameter with each change.  The Last-In-First-Out (LIFO) queue discipline is employed to resolve this issue \cite{Longbo:LIFO}.  However, these works, which achieve average queue size that grows logarithmically in $V$, still assume the availability of infinite buffer space \cite{Neely:EnergyDelay,Longbo:Lagrange,Longbo:LIFO}.

A practical implementation of the LIFO scheme is developed in  \cite{Scott:LIFO}.  The work in \cite{Scott:LIFO} also introduces a \emph{floating-queue} algorithm, operating under the LIFO scheme, to deal with finite buffers.  The algorithm in \cite{Scott:LIFO} is heuristic, and it is not clear how to analyze its behavior.  The current work is inspired by this floating queue idea of \cite{Scott:LIFO} and adopts the same ``floating queue'' terminology, even though the floating-queue algorithm developed here is different from \cite{Scott:LIFO}.  Indeed, the floating queue technique of this paper operates under the First-In-First-Out (FIFO) scheme.  It splits each queue into two queues (one for \emph{real} and one for \emph{fake} packets) and yields analytical guarantees on utility, delay, and packet drops.  

Several backpressure approaches \cite{Ryu:ICNS,Athanasopoulou:randomize,Ji:Delaybased} attempt to improve network delay.  However, those focus on specific aspects and do not have the theoretical utility-delay tradeoff.
Stochastic network optimization with finite buffers has been studied previously in \cite{Long:FiniteBuffer}.  That work uses a non-standard Lyapunov function and knowledge of an $\epsilon$ parameter to derive an upper bound on the required buffer size.  However, the $\epsilon$  parameter can be difficult to determine in practice, and the resulting utility-delay tradeoff is still $[O(1/V), O(V)]$.  An implementation of that work is studied in \cite{Xue:FiniteBufferImplement}.

This paper develops a floating-queue approach to general stochastic network optimization with finite buffers.  Our algorithm is inspired by finite buffer heuristics in \cite{Scott:LIFO} and the  steady state analysis in \cite{Longbo:Lagrange}.  We propose the floating-queue algorithm to solve the learning issue in \cite{Longbo:Lagrange}.  The result obtains the best of both worlds:  It achieves the desired steady state performance but is just as adaptive to network changes as LIFO scheduling. For finite buffers of size $B$, deviation from utility optimality is shown to decrease like $O(e^{-B})$ and packet drops are shown to have rate $O(e^{-B})$, while average per-hop delay is $O(B)$.

This paper is organized as follows.  The system model is described in Section \ref{sec:system_model}. Section \ref{sec:deterministic} describes the standard drift-plus-penalty approach.  The floating queue algorithm is introduced in Section \ref{sec:floatingqueue}.  Performance of the floating-queue algorithm is analyzed in Section \ref{sec:performance} and is validated by simulation in Section \ref{sec:simulation}.  Section \ref{sec:conclusion} concludes the paper.

\section{System Model}
\label{sec:system_model}

The network model of this paper is similar to that of \cite{Longbo:Lagrange}.  Consider a network with $N$ queues that evolve in discrete (slotted) time $t \in \{0, 1, 2, \dotsc \}$.  At each time slot, a \emph{network controller} observes the current \emph{network state} before making a \emph{decision}.  The goal of the controller is to minimize a time average cost subject to network stability.  An example of time average cost is average power incurred over the network.  Utility maximization can be treated by defining the slot-$t$ cost as $-1$ times a slot-$t$ reward. 
An example utility maximization problem is to maximize time average network 
throughput.  The rest of the paper deals with cost minimization, with the understanding that this can also treat utility maximization. 


\subsection{Network State}
The network experiences randomness every time slot.  This randomness is called the \emph{network state} and can 
represent a vector of channel conditions and/or random arrivals for slot $t$. Assume there are $M$ different network states. Define $\set{S} = \{ s_1, s_2, \dotsc, s_M \}$ as the set of all possible states. Let $S(t)$ denote the network state experienced by the network at time $t$.  Let $\pi_m \in [0, 1]$ be the steady state 
probability that $S(t)=s_m$, i.e., $\pi_m = \prob{S(t) = s_m}$.  For simplicity, it is assumed that $S(t)$ is independent and identically distributed (i.i.d.) over slots. The same results can be shown in the general case of ergodic but non-i.i.d. processes (see \cite{Longbo:Lagrange}).
The network controller can observe $S(t)$ before making the slot-$t$ decision, but the $\pi_m$ probabilities are not necessarily known to the controller. 

\subsection{Control Decision}
Every time slot, the network controller chooses a decision from a set of feasible actions which depends on the current network state.  Formally, define $\set{X}^{S(t)}$ as the decision set depending on $S(t)$, and let $x(t)$ denote the decision chosen by the controller at time $t$, where $x(t) \in \set{X}^{S(t)}$.  Assume the action set $\set{X}^{s_m}$ is finite for every $s_m \in \set{S}$.  
On slot $t$, these control decision and network state $(x(t), S(t))$ affects the network in 2 aspects:

1) A cost is incurred.  The cost is $f(t) \defeq f(x(t), S(t)): \set{X}^{S(t)} \rightarrow \RealSet$.  An example cost is energy expenditure. Another example is $-1$ times the amount of newly admitted packets. 

\begin{figure}
  \centering
  \includegraphics[scale=1]{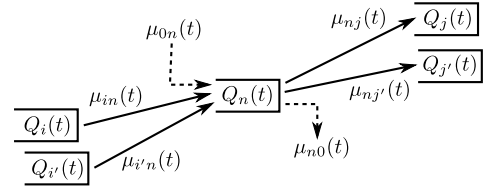}
  \caption{Arrivals and services at a standard queue}
  \label{fig:actionatqueue}
\end{figure}

2) Queues are served. The service variables are $\mu_{ij}(t)$, representing the integer amount of packets taken from queue $i$ and transmitted to queue $j$, for all $i,j \in \mathcal{N} \defeq \{1, \ldots, N\}$.  This is determined by a function $\mu_{ij}(t) \defeq \mu_{ij}(x(t), S(t)): \set{X}^{S(t)} \rightarrow \IntSet_+$.  Further, the decision admits $\mu_{0i}(t) \defeq \mu_{0i}(x(t), S(t)): \set{X}^{S(t)} \rightarrow \IntSet_+$ integer amount of exogenous packets to queue $i \in \set{N}$.  Packets depart from the network at queue $j \in \set{N}$ with an integer amount $\mu_{j0}(t) \defeq \mu_{j0}(x(t), S(t)): \set{X}^{S(t)} \rightarrow \IntSet_+$.  Note that we set $\mu_{ii}(t) = 0$ for all $i \in \set{N}\cup\prtc{0}$ and for all $t$.   The transmission, admission, and departure are shown in Figure \ref{fig:actionatqueue}.

For every $s_m \in \set{S}$, we assume functions $f(\cdot, s_m)$ and $\mu_{ij}(\cdot, s_m)$ for $i,j \in \set{N}\cup\prtc{0}$ are time-invariant, and magnitudes of $\sum_{i=0}^N \mu_{in}(\cdot,s_m)$ and $\sum_{j=0}^N \mu_{nj}(\cdot,s_m)$ are upper bounded by constant $\delta^\maxvar \in (0, \infty)$ for every $n \in \set{N}$.  Furthermore, the network optimization is assumed to satisfy the following \emph{Slater  condition} \cite{Longbo:Lagrange}: Let $\abs{\set{X}}$ be the cardinality of set $\set{X}$.  For every $s_m \in \set{S}$ and $k \in \{1, \dotsc, \abs{\set{X}^{s_m}}\}$, 
there exist probabilities $\zeta_k^{s_m}$ (such that $\sum_{k=1}^{\abs{\set{X}^{s_m}}} \zeta_k^{s_m} = 1$ for all $s_m \in \set{S}$) that define a \emph{stationary and randomized algorithm}.  Whenever the network controller observes state $S(t)=s_m$, the stationary and randomized algorithm chooses action $x_k^{s_m}$ with conditional probability $\zeta_k^{s_m}$. 
The Slater condition assumes there exists such a stationary and randomized algorithm that satisfies: 
\begin{multline*}
  \sum_{m = 1}^M  \sum_{k=1}^{\abs{\set{X}^{s_m}}} \pi_m\zeta_k^{s_m} \prts{ \sum_{i = 0}^N \mu_{in}(x_k^{s_m}, s_m) - \sum_{j=0}^N \mu_{nj}(x_k^{s_m}, s_m) } \\
  \leq -\eta \quad\quad \text{for all}~ n \in \set{N},
\end{multline*}
for some $\eta > 0$. In fact, this assumption is the standard Stater condition of convex optimization \cite{Bertsekas:Convex}.

\subsection{Standard Queue}
\label{sec:standardqueue}
The network consists of $N$ \emph{standard queues}.  Let $Q_n(t)$ denote the backlog in queue $n$ at time $t$, and let $Q(t) \defeq (Q_1(t), \dotsc, Q_N(t))$ be the vector of these backlogs.  The backlog dynamic of queue $n \in \{1, \dotsc, N\}$ is
\begin{equation}
  \label{eq:queue}
  Q_n(t+1) = \max \prts{ Q_n(t) - \sum_{j = 0}^N \mu_{nj}(t), 0 } + \sum_{i = 0}^N \mu_{in}(t).
\end{equation}
When there are not enough packets in a queue, i.e, $Q_n(t) < \sum_{j = 0}^N \mu_{nj}(t)$, blank packets are used to fill up transmissions.

\subsection{Stochastic Formulation}
The controller seeks to minimize the expected time-average cost while maintaing queue stability.
The expected time average cost is defined by
\begin{equation*}
  \bar{f} = \limsup_{T \rightarrow \infty} \frac{1}{T} \sum_{t = 0}^{T-1} \expect{ f(t) },
\end{equation*}
and the queue stability is satisfied when 
\begin{equation*}
  \limsup_{T \rightarrow \infty} \frac{1}{T} \sum_{t = 0}^{T-1} \sum_{n = 1}^N \expect{ Q_n(t) } < \infty.
\end{equation*}

The \emph{stochastic network optimization problem} is
\begin{align}
  \minimize \quad & \bar{f} \label{eq:stoproblem}\\
  \subjectto \quad& \text{queue stability}. \notag
\end{align}


\section{Drift-Plus-Penalty Method}
\label{sec:deterministic}

\subsection{Drift-Plus-Penalty Method}

The drift-plus-penalty method of \cite{Neely:SNO} can solve problem \eqref{eq:stoproblem} via a greedy decision at each time slot that does not require knowledge of the steady state probabilities. The method has parameter $V\geq 0$. In the special case of $V=0$, this policy is also called ``MaxWeight'' or ``Backpressure'':

\textbf{Drift-Plus-Penalty Policy:} At every time $t \in \{0, 1, 2, \dotsc\}$, the network controller observes network state $S(t)$ and backlog vector $Q(t)$.  Decision $x(t) \in \set{X}^{s_m}$ is chosen to solve:
\begin{align}
  \minimize \quad &  V f(x(t), S(t)) \label{eq:dpppolicy}\\
  & \hspace{-3em}  + \sum_{n = 1}^N Q_n(t) \prts{ \sum_{i = 0}^N \mu_{in}(x(t), S(t)) - \sum_{j = 0}^N \mu_{nj}(x(t), S(t)) } \notag\\
  \subjectto \quad & x(t) \in \set{X}^{S(t)}. \notag
\end{align}

Depending on the separability structure of problem \eqref{eq:dpppolicy}, it can be decomposed to smaller subproblems that can be solved distributively.  The algorithm is summarized in Figure \ref{alg:standard}.

\begin{figure}
\begin{boxedalgorithmic}
  \STATE \textbf{Initialization:} $Q(0) = \vect{0}$ \\
  \FOR{$t \in \prtc{0, 1, 2, \dotsc}$}
    \STATE{Observer $S(t)$ and $Q(t)$.} \\
    \STATE{Choose $x(t)$ that solves \eqref{eq:dpppolicy}.} \\
    \STATE{Update $Q_n(t+1)$ according to \eqref{eq:queue} $\forall n \in \set{N}$.} \\
  \ENDFOR
\end{boxedalgorithmic}
  \caption{The drift-plus-penalty algorithm}
  \label{alg:standard}
\end{figure}

It has been shown in \cite{Neely:SNO} that
\begin{align}
  f^\dppvar & \leq f^\optvar + O(1/V) \label{eq:fdpp} \\
  \limsup_{T \rightarrow \infty} \frac{1}{T} \sum_{t = 0}^{T-1} \sum_{n = 1}^N \expect{ Q_n(t) } & = O(V), \label{eq:qdpp}
\end{align}
where $f^\dppvar$ is  the expected time average cost achieved 
by the drift-plus-penalty policy, and $f^\optvar$ is the optimal cost of problem \eqref{eq:stoproblem}.  The inequality \eqref{eq:fdpp} implies that the drift-plus-penalty policy achieves cost within $O(1/V)$ of the optimal cost, which can be made as small as desired by choosing a sufficiently large value of $V$.   The equality \eqref{eq:qdpp} implies that average queue backlog grows linearly with $V$.  Applying Little's law 
gives the $[O(1/V), O(V)]$ utility-delay tradeoff (see \cite{Bertsekas:DataNetworks} for a 
standard description of Little's law). 

Notice that the drift-plus-penalty algorithm assumes infinite buffer size at each queue, even though the 
average queue size is bounded by $O(V)$.

\subsection{Deterministic Problem}
In order to consider a finite buffer regime, the steady-state behavior of the drift-plus-penalty algorithm is considered. 
In \cite{Longbo:Lagrange}, the stochastic problem \eqref{eq:stoproblem} is shown to have an associated deterministic problem as follows: 
\begin{align}
  \minimize \quad & V\sum_{m=1}^M \pi_m f(x^{s_m}, s_m) \label{eq:detproblem}\\
  \subjectto \quad & \sum_{m=1}^M \pi_m \sum_{i=0}^N \mu_{in}(x^{s_m}, s_m) \notag\\
  & \leq \sum_{m=1}^M \pi_m \sum_{j=0}^N \mu_{nj}(x^{s_m}, s_m) ~~ \forall n \in \set{N} \label{eq:detconstraint}\\
  & x^{s_m} \in \set{X}^{s_m} ~~ \forall m \in \{1, \dotsc, M\}. \notag
\end{align}

Let $\dvar = (\dvar_1, \dotsc, \dvar_N)$ be a vector of dual variables associated with constraint \eqref{eq:detconstraint}.
The dual function of problem \eqref{eq:detproblem} is defined as:
\begin{multline}
\label{eq:detdualfunction}
g(\dvar) = \sum_{m=1}^M \pi_m \inf_{x^{s_m} \in \set{X}^{s_m}} \Biggl\{ V f(x^{s_m}, s_m) \\+ \sum_{n=1}^N \dvar_n \Biggl[ \sum_{i=0}^N \mu_{in}(x^{s_m}, s_m) - \sum_{j=0}^N \mu_{nj} (x^{s_m}, s_m) \Biggr]\Biggr\}.
\end{multline}
This dual function \eqref{eq:detdualfunction} is concave.  Therefore, the following 
dual problem is a convex optimization problem:
\begin{align}
  \maximize \quad & g(\dvar) \label{eq:detdualproblem} \\
  \subjectto \quad & \dvar \in \RealSet_+^N. \notag
\end{align}

Let $\dvar^{V\ast} = (\dvar_1^{V\ast}, \dotsc, \dvar_N^{V\ast})$ be a vector of Lagrange multipliers, which solves the dual problem \eqref{eq:detdualproblem} with parameter $V$.  The following theorem from \cite{Longbo:LIFO}
describes a steady state property of the drift-plus-penalty algorithm: 
\begin{theorem}
  \label{thm:steadystate}
  Suppose $\dvar^{V\ast}$ is unique, the Slater condition holds, and the dual function $g(\dvar)$ satisfies:
\begin{equation*}
  g(\dvar^{V\ast}) \geq g(\dvar) + L \norm{ \dvar^{V\ast} - \dvar } \quad \text{for all}~ \dvar \in \RealSet_+^N,
\end{equation*}
for some constant $L > 0$, independent of $V$.  Then under the drift-plus-penalty policy, there exist constants $D, K, c^\ast$, independent of $V$, such that for any $\beta \geq 0$, the following upper bound holds
\begin{equation}
  \label{eq:probbound}
  \set{P}(D, K\beta) \leq c^\ast e^{-\beta},
\end{equation}
where
\begin{multline}
  \label{eq:pdef}
\set{P}(D, K\beta) \\\defeq \limsup_{T \rightarrow \infty} \frac{1}{T} \sum_{t = 0}^{T-1} \prob{ \exists n, \abs{ Q_n(t) - \dvar_n^{V\ast} } > D + K\beta }.
\end{multline}
\end{theorem}
\begin{IEEEproof}
  Please see the full proof in \cite{Longbo:Lagrange}.
\end{IEEEproof}

As all transmissions, admissions, and departures are integers, 
the queue vector has a countably infinite number of 
possibilities,  and under a mild ergodic assumption the steady state distribution of $\prtc{ Q(t) : t \geq 0}$ exists.  In this ergodic case, the  $\set{P}(D,K\beta)$ value in \eqref{eq:pdef} becomes the steady state probability that backlog deviates more than $D+K\beta$ away from the vector of Lagrange multipliers.  Note that 
the probability in \eqref{eq:probbound} vanishes exponentially in  $\beta$. 
This implies that, in steady state, a large portion of arrivals and services occur when the queue backlog vector is close to the Lagrange multiplier vector $\dvar^{V\ast}$.  Thus, if we can admit and serve this portion of traffic using finite-buffer queues, the network still operates near its optimal point.

\section{Floating-Queue Algorithm}
\label{sec:floatingqueue}
In this section, the floating-queue algorithm is presented as a way to implement the drift-plus-penalty algorithm using
finite buffers.  The algorithm preserves the dynamics of the drift-plus-penalty algorithm and hence inherits several of its performance guarantees.

Recall that standard queue $n \in \set{N}$ has dynamic \eqref{eq:queue}.  To simplify notation, let $a_n(t) \defeq \sum_{i=0}^N \mu_{in}(t)$ denote aggregated arrivals to queue $n$, and let $b_n(t) \defeq \sum_{j=0}^N \mu_{nj}(t)$ denote aggregated services from queue $n$ at time $t$.  This implies that $\delta^\maxvar$ upper bounds both $a_n(t)$ and $b_n(t)$.  The dynamic \eqref{eq:queue} can be written as
\begin{equation}
  \label{eq:queue_simple}
  Q_n(t+1) = \max\prts{ Q_n(t) - b_n(t), 0 } + a_n(t).
\end{equation}
For the rest of this paper, the above dynamic is considered for a standard queue.  Note that $a_n(t)$ and $b_n(t)$ are fully determined after knowing all $\mu_{ij}(t)$ from the drift-plus-penalty algorithm.

\subsection{Queue Transformation}
In the floating-queue algorithm, each standard queue $n \in \set{N}$ of Section \ref{sec:standardqueue} is a combination of a \emph{real queue} and a \emph{fake queue}.  The real queue has buffer size $B$ for storing \emph{real packets}.  The fake queue contains \emph{fake packets} and only requires a counter to implement.  
Let $Q_n^r(t)$ and $Q_n^f(t)$ denote respectively the amount of backlogs in the real queue and the fake queue of the standard queue $n$.  Define $Q^r(t) \defeq (Q_1^r(t), \dotsc, Q_N^r(t))$ and $Q^f(t) \defeq (Q_1^f(t), \dotsc, Q_N^f(t))$ as vectors of real and fake queue backlogs.  We use the term \emph{floating queue $n$} to refer to the 2-queue combination consisting of real and fake queues $n$.


\begin{figure}
  \centering
  \includegraphics[scale=0.9]{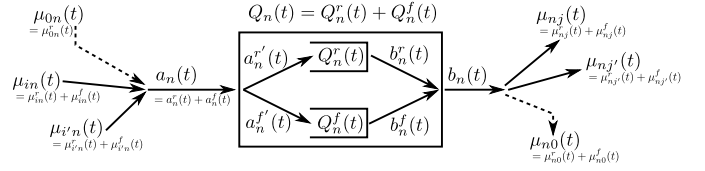}
  \caption{Transformation of a standard queue to a floating queue}
  \label{fig:floatingqueue}
\end{figure}

\subsection{Real and Fake Parts of Arrivals and Services}
At each queue $n \in \set{N}$, let $b_n^r(t)$ and $b_n^f(t)$ denote the aggregated real and fake serviced packets at time $t$.  the floating-queue algorithm always serves real packets before fake packets as:
\begin{align}
  b_n^r(t) &= \min\prts{ Q_n^r(t), b_n(t) } \label{eq:bnr} \\
  b_n^f(t) &= b_n(t) - b_n^r(t). \label{eq:bnf}
\end{align}
It is easy to see that $b_n(t) = b_n^r(t) + b_n^f(t)$.  Also, $b_n^r(t)$ and $b_n^f(t)$ are fully determined, since $b_n(t)$ and $Q_n^r(t)$ are known.

To differentiate the real and fake packets in the drift-plus-penalty variable $\mu_{ij}(t)$, let $\mu_{ij}^r(t)$ and $\mu_{ij}^f(t)$ denote the real and fake parts of $\mu_{ij}(t)$ for $i,j \in \set{N}\cup\prtc{0}$.  All $\mu_{ij}(t)$ are non-negative integers.  Since $\mu_{ii}(t) = 0$, we have $\mu_{ii}^r(t) = \mu_{ii}^f(t) = 0$ for all $i \in \set{N}\cup\prtc{0}$.  All exogenous arrivals are considered as real packets, so knowing the value of $\mu_{0j}(t)$, we set $\mu_{0j}^r(t) = \mu_{0j}(t)$ and $\mu_{0j}^f(t) = 0$ for every $j \in \{1, \ldots, N\}$.  The real and fake parts of $\mu_{ij}(t)$ for $i \in \set{N}, j \in \set{N}\cup\{0\}$ can be set arbitrarily to satisfy:
\begin{align*}
b_n^r(t) = \sum_{j=0}^N \mu_{nj}^r(t), \quad b_n^f(t) = \sum_{j=0}^N \mu_{nj}^f(t) \quad\quad\forall n \in \set{N} \\
\mu_{ij}(t) = \mu_{ij}^r(t) + \mu_{ij}^f(t) \quad\quad \forall i \in \set{N}, j \in \set{N}\cup\prtc{0}
\end{align*}
Therefore, all $\mu_{ij}^r(t)$ and $\mu_{ij}^f(t)$ for all $i,j \in \set{N}\cup\prtc{0}$ are fully determined.  This is illustrated in Figure \ref{fig:floatingqueue}.

Let $a_n^r(t) = \sum_{i=0}^N \mu_{in}^r(t)$ and $a_n^f(t) = \sum_{i=1}^N \mu_{in}^f(t)$ be the aggregated real and fake parts of $a_n(t)$.  They are fully determined and $a_n(t) = a_n^r(t) + a_n^f(t)$.

The floating-queue algorithm always admits a real packet to a real queue as much as allowed by its buffer space.  Let $a_n^{r'}(t)$ denote the amount of packets in $a_n^r(t)$ that are admitted to the real queue $n$ at time $t$.  Another part of $a_n^r(t)$, which is dropped, is denoted by $d_n(t)$ and becomes fake packets.  Let $a_n^{f'}(t)$ denote the admitted fake arrivals, including original fake arrivals $a_n^f(t)$ and dropped packets $d_n(t)$.  The arrival dynamic of queue $n \in \set{N}$ is
\begin{align}
  a_n^{r'}(t) &= \min\prts{B - Q_n^r(t), a_n^r(t)} \label{eq:anr} \\
  d_n(t)     &= a_n^r(t) - a_n^{r'}(t) \label{eq:dn} \\
  a_n^{f'}(t) &= a_n^f(t) + d_n(t). \label{eq:anf}
\end{align}
It is easy to see that $a_n^{r'}(t) + a_n^{f'}(t) = a_n^r(t) + a_n^f(t)$.

\subsection{Real and Fake Queuing Dynamics}
The dynamics of real and fake queues of standard queue $n \in \set{N}$ are
\begin{align}
  Q_n^r(t+1) &= Q_n^r(t) - b_n^{r}(t) + a_n^{r'}(t) \label{eq:qnr}\\
  Q_n^f(t+1) &= \max\prts{ Q_n^f(t) - b_n^f(t), 0 } + a_n^{f'}(t). \label{eq:qnf}
\end{align}

\begin{lemma}
  \label{lem:qsync}
  From any time $t_0$, when $Q(t_0) = Q^r(t_0) + Q^f(t_0)$, it follows that
\begin{equation*}
  Q(t) = Q^r(t) + Q^f(t), \quad \forall t \geq t_0.
\end{equation*}
\end{lemma}
\begin{IEEEproof}
We prove this lemma by induction. For $t_0$, $Q(t_0) = Q^r(t_0) + Q^f(t_0)$ by the assumption.  Suppose $Q(t) = Q^r(t) + Q^f(t)$ at time $t$.  At queue $n \in \set{N}$, we have
\begin{align}
  Q_n(t+1) & = \max\prts{ Q_n(t) - b_n(t), 0 } + a_n(t) \notag\\
  & = \max\prts{ Q_n^r(t) + Q_n^f(t) - b_n^{r}(t) - b_n^{f}(t), 0 } \notag\\
  & \quad+ a_n^{r'}(t) + a_n^{f'}(t). \label{eq:la1}
\end{align}

When there are not enough real packets, $Q_n^r(t) < b_n(t)$, it follows that $Q_n^r(t) - b_n^r(t) = 0$ from \eqref{eq:bnr}. Equation \eqref{eq:la1} becomes
\begin{align*}
  Q_n(t+1) & = \max\prts{ Q_n^f(t) - b_n^f(t), 0 } + a_n^{f'}(t) + a_n^{r'}(t) + 0 \\
  & = Q_n^f(t+1) + Q_n^r(t+1).
\end{align*}

When there are enough real packets, $Q_n^r(t) \geq b_n(t)$, we have that $b_n^f(t) = 0$ from \eqref{eq:bnr} and \eqref{eq:bnf}. Equation \eqref{eq:la1} becomes
\begin{align*}
  Q_n(t+1) & = Q_n^f(t) + a_n^{f'}(t) + Q_n^r(t) - b_n^r(t) + a_n^{r'}(t) \\
  & = Q_n^f(t+1) + Q_n^r(t+1).
\end{align*}
Thus, $Q_n(t+1) = Q_n^f(t) + Q_n^r(t)$ for all $n \in \set{N}$.
\end{IEEEproof}

The implication of Lemma \ref{lem:qsync} is that, although the floating-queue algorithm implements these real and fake queues instead of the standard queues, the dynamics of $Q(t)$ and $Q^r(t) + Q^f(t)$ are the same.  Hence, when decision $x(t)$ is chosen by solving \eqref{eq:dpppolicy} with $Q^r(t) + Q^f(t)$ instead of $Q(t)$, all decisions $\prtc{ x(t) }_{t=0}^\infty$ under the standard algorithm (in Figure \ref{alg:standard}) are identical to the decisions $\prtc{ x(t) }_{t=0}^\infty$ under the floating-queue algorithm (in Figure \ref{alg:FQ}), given that $Q^r(0) + Q^f(0) = Q(0)$.  Yet, the buffer size of each real queue in the floating-queue algorithm is $B$.  Let $Q^f(0) \in \IntSet_+^N$ be a predefined fake backlogs.  The floating-queue algorithm is summarized in Figure \ref{alg:FQ}.

\begin{figure}
\begin{boxedalgorithmic}
  \STATE \textbf{Initialization:} $Q^r(0) = \vect{0}$ and $Q^f(0)$ \\
  \FOR{$t \in \prtc{0, 1, 2, \dotsc}$}
    \STATE{Observer $S(t)$ and let $Q(t) = Q^r(t) + Q^f(t)$.} \\
    \STATE{Choose $x(t)$ that solves \eqref{eq:dpppolicy}.} \\
    \STATE{Calculate $(a_n(t), b_n(t)) ~ \forall n \in \set{N}$.}
    \STATE{Calculate $(b_n^r(t), b_n^f(t))$ as \eqref{eq:bnr} and \eqref{eq:bnf} $\forall n \in \set{N}$.} \\
    \STATE{Adjust $(a_n^{r'}(t), a_n^{f'}(t)) $ as \eqref{eq:anr}--\eqref{eq:anf} $\forall n \in \set{N}$.} \\
    \STATE{$Q_n^r(t+1) = Q_n^r(t) - b_n^{r'}(t) + a_n^{r'}(t) ~ \forall n \in \set{N}$.} \\
    \STATE{$Q_n^f(t+1) = \max\prts{ Q_n^f(t) - b_n^f(t), 0 } + a_n^{f'}(t) ~ \forall n \in \set{N}$.} \\
  \ENDFOR
\end{boxedalgorithmic}
  \caption{The floating-queue algorithm}
  \label{alg:FQ}
\end{figure}

We prove a useful lemma of the floating-queue algorithm, which will be used in Section \ref{sec:samplepathanalysis}.
\begin{lemma}
  \label{lem:dropincrease}
  Under the floating-queue algorithm, when the buffer size of the real queue $n \in \set{N}$ is $B \geq 2\delta^\maxvar$, if $d_n(t) > 0$, then $Q_n^f(t+1) > Q_n^f(t)$.
\end{lemma}
\begin{IEEEproof}
  Event $d_n(t) > 0$ implies that $a_n^r(t) > a_n^{r'}(t)$ from \eqref{eq:dn} and $a_n^{r'}(t) = B - Q_n^r(t)$ from \eqref{eq:anr}, so $Q_n^r(t) > B - a_n^r(t)$.  When $B \geq 2\delta^\maxvar$, we have $Q_n^r(t) > 2\delta^\maxvar - a_n^r(t) \geq \delta^\maxvar$, and there are enough real packets for all services.  Therefore, all services take real packets and $b_n^r(t) = b_n(t)$ and $b_n^f(t) = 0$ from \eqref{eq:bnr} and \eqref{eq:bnf}.  From \eqref{eq:qnf} and \eqref{eq:anf}, we have 
\begin{equation*}
  Q_n^f(t+1) = Q_n^f(t) + a_n^{f'}(t) = Q_n^f(t) + a_n^f(t) + d_n(t) > Q_n^f(t).
\end{equation*}
\end{IEEEproof}
The interpretation of Lemma \ref{lem:dropincrease} is that, for any queue $n$ with buffer size $B \geq 2\delta^\maxvar$, if real packets are dropped at time $t$, then the fake backlogs at time $t+1$ always increase.

\section{Performance Analysis}
\label{sec:performance}
The steady-state performance of the floating-queue algorithm is analyzed by bounding from below the admitted real arrivals at each queue $n \in \set{N}$.  Define $\spath_n(t) \defeq$
\begin{equation*}
  \prtr{ a_n^r(t), a_n^f(t), a_n^{r'}(t), a_n^{f'}(t), b_n^r(t), b_n^f(t), Q_n^r(t), Q_n^f(t) }
\end{equation*}
as a sample path of the arrivals, services, and backlogs of queue $n$ that is generated by the floating-queue algorithm at time $t$.  For any positive integer $T$ and starting time $t_0$, a sample path of queue $n$ from $t_0$ to $t_0+T$ is denoted by $\prtc{ \spath_n(t) }_{t = t_0}^{t_0+T}$.  Note that $Q_n(t)$ can be recovered from this sample path as $Q_n(t) = Q_n^r(t) + Q_n^f(t)$.

From sample path $\spath_n(t)$, the amount of real arrivals are $a_n^r(t)$, and the amount of \emph{admitted} real arrivals are $a_n^{r'}(t)$, which depend on the floating-queue mechanism \eqref{eq:anr}.  To lower bound this admitted real arrival $a_n^{r'}(t)$, we construct another mechanism, called ``lower-bound policy'', that operates over the sample path.  It has a different rule for counting admitted real packets (later defined as $\hat{a}_n^r(t)$), which is part of the real arrivals $a_n^r(t)$.  We will show (Lemma \ref{lem:THTLgeneral}) that the amount of admitted real arrivals under the floating-queue algorithm is lower bounded by the amount of admitted real arrivals under the lower-bound policy.  Using this lower bound the performance of the floating-queue algorithm can be analyzed.

\subsection{Lower-Bound Policy}
\label{sec:lowerboundpolicy}
In this section, queue $n \in \set{N}$ is fixed and the lower-bound policy is defined for this queue.  For simplicity, assume that the buffer size $B$ is even and $B \geq 2\delta^\maxvar$.  

Recall that $\dvar^{V\ast}$ is the Lagrange multiplier of problem \eqref{eq:detdualproblem}, and $\delta^\maxvar$ is the upper bound on $a_n(t)$ and $b_n(t)$.  Define
\begin{equation*}
  \set{B}_n \defeq \prts{ \dvar_n^{V\ast} - B/2 + \delta^\maxvar, \dvar_n^{V\ast} + B/2 - \delta^\maxvar }.
\end{equation*}

Let $\hat{a}_n^r(t)$ denote the number of admitted real packets under the lower-bound policy at time $t$.  Given any sample path $\spath_n(t)$ having real arrivals $a_n^r(t)$ and total backlogs $Q_n(t)$, the lower-bound policy counts real packets as
\begin{equation}
  \label{eq:hanr}
  \hat{a}_n^r(t) = \left\{
    \begin{array}{ll}
      a_n^r(t) & , Q_n(t) \in \set{B}_n \\
      0  & , Q_n(t) \notin \set{B}_n.
    \end{array}\right.
\end{equation}

Let $\hat{d}_n(t)$ denote the number of dropped packets under the lower-bound policy at time $t$.  It satisfies
\begin{equation}
  \label{eq:hdn}
  \hat{d}_n(t) = a_n^r(t) - \hat{a}_n^{r}(t).
\end{equation}
Notice that $\hat{a}_n^r(t)$ and $\hat{d}_n(t)$ are artificial numbers and are not real and fake packets in a real system.  These values can be determined by $a_n^r(t)$ and $Q_n(t)$ of sample path $\spath_n(t)$.

\subsection{Sample Path Analysis}
\label{sec:samplepathanalysis}
The goal of this section is to show (Lemma \ref{lem:THTLgeneral}) that, for any sample path $\prtc{ \spath_n(t) }_{t=t_0}^\infty$ of queue $n \in \set{N}$ and any positive integer $T$, the admitted real arrivals under the floating-queue algorithm with buffer size $B$ is lower bounded by
\begin{equation*}
  \sum_{t = t_0}^{t_0+T-1} a_n^{r'}(t) \geq \sum_{t = t_0}^{t_0+T-1} \hat{a}_n^{r}(t) - B.
\end{equation*}

In this section, queue $n$ is fixed and analyzed; however, the analysis results hold for every queue $n \in \set{N}$.

\begin{figure}
  \centering
  \includegraphics[scale=0.95]{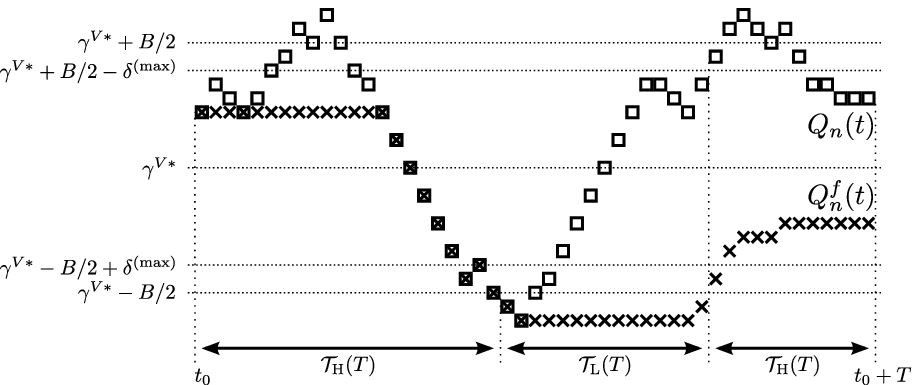}
  \vspace{-2em}
  \caption{Time interval $\set{T}(T)$ is partitioned into $\set{T}_\HH(T)$ and $\set{T}_\LL(T)$.}
  \label{fig:THTL}
\end{figure}

For any starting time $t_0$ and any positive integer $T$, define $\set{T}(T) \defeq \prtc{t_0, \dotsc, t_0+T}$ as a time interval of consideration.  It can be partitioned into disjoint sets $\set{T}_\HH(T)$ and $\set{T}_\LL(T)$, which are illustrated in Figure \ref{fig:THTL}, where
\begin{align*}
  \set{T}_\HH(T) &\defeq \prtc{ t \in \set{T}(T) : Q_n^f(t) \geq \dvar_n^{V\ast} - B/2 } \\
  \set{T}_\LL(T) &\defeq \prtc{ t \in \set{T}(T) : Q_n^f(t) < \dvar_n^{V\ast} - B/2 }.
\end{align*}

\begin{figure}
  \centering
  \includegraphics[scale=0.95]{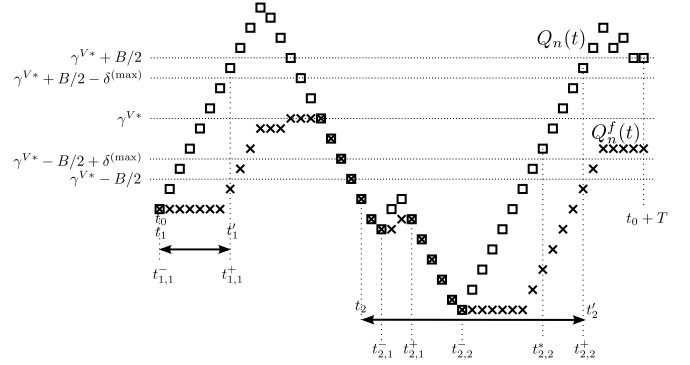}
  \vspace{-0.5em}
  \caption{Partitioning $\set{T}_\LL(T)$ into intervals of $t_k$, $t_k'$, $t_{k,j}^-$, and $t_{k,j}^+$}
  \label{fig:TL}
\end{figure}

Time interval $\set{T}_\LL(T)$ can be partitioned into disjointed intervals of time that starts when the fake queue in the sample path satisfies $Q_n^f(t) < \dvar_n^{V\ast} - B/2$ and ends when it does not.  This is illustrated in Figure \ref{fig:TL}.  For $k \in \{1, 2, \dotsc\}$, let
\begin{align*}
  t_k &= \arginf_{t \in \{t_{k-1}' +1, \dotsc, t_0+T\}} \prtc{ Q_n^f(t) < \dvar_n^{V\ast} - B/2 } \\
  t_k' &= \arginf_{t \in \{t_k +1, \dotsc, t_0+T\}} \prtc{ Q_n^f(t) \geq \dvar_n^{V\ast} - B/2 } - 1,
\end{align*}
where $t_0' = t_0-1$ and $\arginf_{t \in \{A, \dotsc, B\}} \{ C(t) \} = B+1$ if $A > B$ or $C(t)$ is not satisfied for all $t \in \{A, \dotsc, B\}$.  Let $K(T) = \argmax_{k \geq 0} \prtc{ t_k < t_0 +T + 1 }$ denote the number of intervals $\{t_k, \dotsc, t_k'\}$ contained in $\set{T}(T)$.

When $K(T) > 0$, the time interval in interval $\{t_k, \dotsc, t_k'\}$ for $k \in \{1, 2, \dotsc, K(T)\}$ can be partitioned into intervals between local minima and local maxima.  Define $U(t) \defeq \arginf_{\tau \in \{t+1, t+2, \dotsc\}} \prtc{ Q_n^f(\tau) > Q_n^f(t) }$ to be the first time index after $t$ that the fake queue increases.  For $k \in \{1, 2, \dotsc, K(T) \}, j \in \{1, 2, \dotsc\}$, let
\begin{align*}
  t_{k,j}^- &= \min\biggl[ \arginf_{t \in \{t_{k,j-1}^+ + 1, \dotsc, t_k' \}} \bigl\{ Q_n^f(t) < Q_n^f(t-1) \text{ and }\\
    & \quad\quad\quad  Q_n^f(t) \leq Q_n^f(\tau)~ \forall \tau \in \{ t+1, \dotsc, U(t) \} \bigr\} , t_k' \biggr] \\
  t_{k,j}^+ &= \min\prts{ \arginf_{t \in \{t_{k,j}^- + 1, \dotsc, t_k' \}} \prtc{ Q_n^f(t) > Q_n^f(t+1) } , t_k' },
\end{align*}
where $t_{k,0}^+ = t_k -1$ and $Q_n^f(t_0 - 1) = \infty$.  Intuitively, during $\{t_k, \dotsc, t_k'\}$, $t_{k,j}^-$ is the first time index that the $j^\text{th}$ local minimum is reached, and $t_{k,j}^+$ is the last time index of the $j^\text{th}$ local maximum.  This is illustrated in Figure \ref{fig:TL}.  Let $J(k) = \arginf_{j > 0} \prtc{ t_{k,j}^+ = t_k' }$ denote the number of local maxima during $\{ t_k, \dotsc, t_k' \}$.

For a technical reason (used in Lemma \ref{lem:TL}), let $\set{T}_\KK(T) \defeq \{(t_k-1) \in \set{T}(T) : k \in \{1, \dotsc, K(T) \}$.  The following lemmas hold for the real arrivals in $\set{T}_\HH(T)\backslash\set{T}_\KK(T)$ and $\set{T}_\LL(T)\cup\set{T}_\KK(T)$.
\begin{lemma}
  \label{lem:TH}
  When $B \geq 2\delta^\maxvar$, given any sample path $\prtc{ \spath_n(t) }_{t = t_0}^{t_0 + T}$, the following relation holds
\begin{equation*}
  \sum_{t \in \set{T}_\HH(T)\backslash \set{T}_\KK(T)} a_n^{r'}(t) \geq \sum_{t \in \set{T}_\HH(T)\backslash\set{T}_\KK(T)} \hat{a}_n^{r}(t).
\end{equation*}
\end{lemma}
\begin{IEEEproof}
Two cases are examined.

1) When $Q_n(t) \in \set{B}_n$ for any $t \in \set{T}_\HH(T)\backslash \set{T}_\KK(T)$, we have $a_n^{r'}(t) = a_n^r(t)$, because real queue $n$ has enough buffer space:
\begin{align*}
  Q_n^r(t) &= Q_n(t) - Q_n^f(t) \\
  &\leq \prtr{ \dvar_n^{V\ast} + B/2 - \delta^\maxvar } - \prtr{ \dvar_n^{V\ast} - B/2 } \\
  &\leq B - \delta^\maxvar.
\end{align*}
The first inequality holds because of $Q_n(t) \in \set{B}_n$ and $t \in \set{T}_\HH(T)$.
For the lower-bound policy, we have $\hat{a}_n^{r}(t) = a_n^r(t)$, because $Q_n(t) \in \set{B}_n$.  So $a_n^{r'}(t) = \hat{a}_n^r(t)$.

2) When $Q_n(t) \notin \set{B}_n$ for any $t \in \set{T}_\HH(T)\backslash \set{T}_\KK(T)$, we have $a_n^{r'}(t) \geq \hat{a}_n^{r}(t) = 0$, because $Q_n(t) \notin \set{B}_n$ and $a_n^{r'}(t) \geq 0$.

These two cases implies the lemma.
\end{IEEEproof}

\begin{lemma}
  \label{lem:TL}
  When $B \geq 2\delta^\maxvar$, given sample path $\prtc{ \spath_n(t) }_{t = t_0}^{t_0 + T}$ with $Q_n^r(t_0) = 0$, the following holds
\begin{equation*}
  \sum_{t \in \set{T}_\LL(T)\cup\set{T}_\KK(T)} a_n^{r'}(t) \geq \sum_{t \in \set{T}_\LL(T)\cup\set{T}_\KK(T)} \hat{a}_n^{r}(t).
\end{equation*}
\end{lemma}
\begin{IEEEproof}
For $k \in \{1, \dotsc, K(T)\}, j \in \{1, \dotsc, J(k)\}$, 3 cases are examined.

1) For real arrivals $a_n^r(t)$ during $t \in \prtc{ t_k -1, \dotsc, t_{k,1}^- -2 }$ (if exists), the fake backlogs $Q_n^f(t)$ is non-increasing by the definition of $t_{k,1}^-$.  From Lemma \ref{lem:dropincrease}, the non-increasing implies no packet drops and $a_n^{r'}(t) = a_n^r(t)$ for $t \in \prtc{ t_k -1, \dotsc, t_{k,1}^- -2 }$.  Since $\hat{a}_n^r(t) \leq a_n^r(t)$, it follows that $\sum_{t=t_k -1}^{t_{k,1}^- -2} a_n^{r'}(t) \geq \sum_{t=t_k -1}^{t_{k,1}^- -2} \hat{a}_n^{r}(t)$.  For a special case when $t_1 = t_0$, same argument can be used to obtain $\sum_{t=t_1}^{t_{1,1}^- -2} a_n^{r'}(t) \geq \sum_{t=t_1}^{t_{1,1}^- -2} \hat{a}_n^{r}(t)$.

2) For real arrivals $a_n^r(t)$ during $t \in \prtc{ t_{k,j}^- -1, \dotsc, t_{k,j}^+ }$, it can be shown that $\sum_{t=t_{k,j}^- -1}^{t_{k,j}^+} a_n^{r'}(t) \geq \sum_{t=t_{k,j}^- -1}^{t_{k,j}^+} \hat{a}_n^{r}(t)$.  This result is proven in Appendix.

3) For real arrivals $a_n^r(t)$ during $t \in \prtc{ t_{k,j}^+ + 1, \dotsc, t_{k,j+1}^- -2}$ (if exists), the fake backlogs $Q_n^f(t)$ is non-increasing by the definitions of $t_{k,j}^+$ and $t_{k,j+1}^-$.  Lemma \ref{lem:dropincrease} implies that $\sum_{t=t_{k,j}^+ +1}^{t_{k,j+1}^- -2} a_n^{r'}(t) \geq \sum_{t=t_{k,j}^+ +1}^{t_{k,j+1}^- -2} \hat{a}_n^{r}(t)$.

Since $Q_n^r(t_0) = 0$, the time interval $\prtc{t_k, \dotsc, t_k'}$ for $k \in \prtc{1, \dotsc, K(T)}$ starts with either the first case or the second case and ends with the second case.  Thus, it is clear that, if real arrivals under the lower-bound policy are dominated over every subinterval, then they are also dominated over the union of these subintervals.
\end{IEEEproof}

\begin{lemma}
  \label{lem:THTL}
  When $B \geq 2\delta^\maxvar$, given sample path $\prtc{ \spath_n(t) }_{t = t_0}^\infty$ with $Q_n^r(t_0) = 0$ and positive integer $T$, it holds that
\begin{equation*}
  \sum_{t = t_0}^{t_0+T-1} a_n^{r'}(t) \geq \sum_{t = t_0}^{t_0+T-1} \hat{a}_n^{r}(t).
\end{equation*}
\end{lemma}
\begin{IEEEproof}
Disjoint time intervals $\set{T}_\HH(T-1)\backslash\set{T}_\KK(T-1)$ and $\set{T}_\LL(T-1)\cup\set{T}_\KK(T-1)$ are the partitions of $\set{T}(T-1)$.  Then Lemma \ref{lem:TH} and Lemma \ref{lem:TL} imply the lemma.
\end{IEEEproof}

The above lemma is not general, since it requires $Q_n^r(t_0) = 0$.  Now its general version is provided.
\begin{lemma}
  \label{lem:THTLgeneral}
  When $B \geq 2\delta^\maxvar$, given any sample paths $\prtc{ \spath_n(t) }_{t=t_0}^\infty$ and positive integer $T$, it holds for any $Q_n^r(t_0) \in \prtc{0, 1, \dotsc, B}$ that
\begin{equation*}
  \sum_{t = t_0}^{t_0+T-1} a_n^{r'}(t) \geq \sum_{t = t_0}^{t_0+T-1} \hat{a}_n^{r}(t) - B.
\end{equation*}
\end{lemma}
\begin{IEEEproof}
Construct sample path $\prtc{\tilde{\spath}_n(t)}_{t_{-1}}^\infty$ with $\tilde{\spath}_n(t) = \prtr{ \tilde{a}_n^r(t), \tilde{a}_n^f(t), \tilde{a}_n^{r'}(t), \tilde{a}_n^{f'}(t), \tilde{b}_n^r(t), \tilde{b}_n^f(t), \tilde{Q}_n^r(t), \tilde{Q}_n^f(t) }$ and
\begin{itemize}
\item $t_{-1} < t_0$,
\item $\tilde{\spath}_n(t) = \spath_n(t)$ for all $t \in \{ t_0, t_0+1, \dotsc \}$,
\item $\tilde{Q}_n^f(t) = Q_n^f(t_0)$ for all $t \in \{t_{-1}, \dotsc, t_0-1\}$,
\item $\tilde{Q}_n^r(t_{-1}) = 0$, $\tilde{Q}_n^r(t_0) = Q_n^r(t_0) = \sum_{t = t_{-1}}^{t_0 -1} \tilde{a}_n^{r'}(t)$,
\item $\tilde{a}_n^{f'}(t) = \tilde{b}_n^r(t) = \tilde{b}_n^f(t) = 0$ for all $t \in \{t_{-1}, \dotsc, t_0-1 \}$.
\end{itemize}
The last two conditions automatically set the values of $\prtc{\tilde{a}_n^{r}(t), \tilde{a}_n^{f}(t), \tilde{Q}_n^r(t) }_{t_{-1}}^{t_0-1}$.  This new sample path satisfies Lemma \ref{lem:THTL}.  Let $\hat{\tilde{a}}_n^{r}(t)$ denote the admitted real arrivals under the lower-bound policy of the new sample path.  Since $\sum_{t = t_0}^{t_0 + T -1} \tilde{a}_n^{r'}(t) = \sum_{t = t_{0}}^{t_0 + T -1} a_n^{r'}(t)$, it follows that
\begin{align*}
  & \sum_{t = t_{-1}}^{t_0-1} \tilde{a}_n^{r'}(t) + \sum_{t = t_{0}}^{t_0 + T -1} a_n^{r'}(t) = \sum_{t = t_{-1}}^{t_0-1} \tilde{a}_n^{r'}(t) + \sum_{t = t_{0}}^{t_0 + T -1} \tilde{a}_n^{r'}(t) \\
  & \geq \sum_{t = t_{-1}}^{t_0-1} \hat{\tilde{a}}_n^{r}(t) + \sum_{t = t_{0}}^{t_0 + T -1} \hat{\tilde{a}}_n^{r}(t) = \sum_{t = t_{-1}}^{t_0-1} \hat{\tilde{a}}_n^{r}(t) + \sum_{t = t_{0}}^{t_0 + T -1} \hat{a}_n^{r}(t).
\end{align*}
The first inequality is the application of Lemma \ref{lem:THTL}, and the last equality holds, because $\prtc{ Q_n(t) }_{t=t_0}^\infty$ of both original and new sample paths are identical.

Therefore, we have
\begin{equation*}
  \sum_{t = t_{0}}^{t_0 + T -1} a_n^{r'}(t) \geq \sum_{t = t_{0}}^{t_0 + T -1} \hat{a}_n^{r}(t) - B.
\end{equation*}
The inequality uses the facts that $\sum_{t = t_{-1}}^{t_0-1} \hat{\tilde{a}}_n^{r}(t) \geq 0$ and $\sum_{t = t_{-1}}^{t_0-1} \tilde{a}_n^{r'}(t) = Q_n^r(t_0) \leq B$.
\end{IEEEproof}

\subsection{Performance of Floating-Queue Algorithm}

\subsubsection{Average Drops}
The average drops at each queue is analyzed using the steady state and sample path results.
Recall that constants $D, K, c^\ast$ are defined in Theorem \ref{thm:steadystate}.

\begin{lemma}
  \label{lem:drop}
  Suppose $B > 2(\delta^\maxvar+D)$.  In the steady state, the average drops at real queue $n \in \set{N}$ under the floating-queue algorithm is bounded by
\begin{equation*}
  \lim_{T \rightarrow \infty} \frac{1}{T} \sum_{t = t_0}^{t_0+T-1} \expect{d_n(t)} \leq \delta^\maxvar c^\ast e^\frac{-[B/2 - \delta^\maxvar -D]}{K}.
\end{equation*}
\end{lemma}
\begin{IEEEproof}
We consider queue $n \in \set{N}$.  Let $\indicator{A}$ be an indicator function of statement $A$ such that $\indicator{A} = 1$ if statement $A$ is true; otherwise $\indicator{A} = 0$.  Equation \eqref{eq:hanr} can be written as $\hat{a}_n^{r}(t) = a_n^r(t) \indicator{ Q_n(t) \in \set{B}_n } = a_n^{r}(t) - a_n^r(t) \indicator{ Q_n(t) \notin \set{B}_n }$.  Then we have that
\begin{multline*}
  \frac{1}{T} \sum_{t=t_0}^{t_0+T-1} \expect{ \hat{a}_n^{r}(t) } \\
  = \frac{1}{T} \sum_{t=t_0}^{t_0+T-1} \expect{ a_n^{r}(t) - a_n^{r}(t)  \indicator{ Q_n(t) \notin \set{B}_n } }
\end{multline*}

Dividing the result in Lemma \ref{lem:THTLgeneral} by $T$ and taking an expectation yields
\begin{equation*}
  \frac{1}{T} \sum_{t = t_0}^{t_0+T-1} \expect{a_n^{r'}(t)} \geq \frac{1}{T} \sum_{t = t_0}^{t_0+T-1} \expect{\hat{a}_n^{r}(t)} - \frac{B}{T}.
\end{equation*}
Combining the above two relations gives
\begin{multline*}
  \frac{1}{T} \sum_{t = t_0}^{t_0+T-1} \expect{a_n^{r'}(t)} \geq \frac{1}{T} \sum_{t=t_0}^{t_0+T-1} \expect{ a_n^{r}(t) }  - \frac{B}{T} \\
  - \frac{1}{T} \sum_{t = t_0}^{t_0+T-1} \expect{a_n^{r}(t)  \indicator{ Q_n(t) \notin \set{B}_n } }.
\end{multline*}
It follow from \eqref{eq:dn} that
\begin{align*}
  \frac{1}{T} \sum_{t = t_0}^{t_0+T-1} \expect{d_n(t)} 
  \leq \frac{1}{T} \sum_{t = t_0}^{t_0+T-1} \delta^\maxvar \prob{ Q_n(t) \notin \set{B}_n } + \frac{B}{T}
\end{align*}
Taking limit as $T$ approaches infinity yields
\begin{multline}
  \label{eq:lb1}
  \lim_{T \rightarrow \infty} \frac{1}{T} \sum_{t = t_0}^{t_0+T-1} \expect{d_n(t)} \\\leq \delta^\maxvar \lim_{T \rightarrow \infty} \frac{1}{T} \sum_{t = t_0}^{t_0+T-1} \prob{ Q_n(t) \notin \set{B}_n }.
\end{multline}
In steady state, Theorem \ref{thm:steadystate} with $\beta = \frac{B/2 - \delta^\maxvar -D}{K}$ yields
\begin{align*}
  & \lim_{T \rightarrow \infty} \frac{1}{T} \sum_{t=t_0}^{t_0+T-1} \prob{  Q_n(t) \notin \set{B}_n } \notag\\
  & \hspace{-0.3em}\leq \limsup_{T \rightarrow \infty} \frac{1}{T} \sum_{t = t_0}^{t_0+T-1} \prob{ \exists n, \abs{ Q_n(t) - \dvar_n^{V\ast} } > B/2 - \delta^\maxvar } \notag\\
  & \hspace{-0.3em}= \set{P}(D, B/2 -\delta^\maxvar - D) \leq c^\ast e^{-[B/2 -\delta^\maxvar -D]/K}.
\end{align*}
Applying the above bound to \eqref{eq:lb1} proves the lemma.
\end{IEEEproof}

\subsubsection{Delay}
At each queue, the average delay experienced by real packets is derived by invoking Little's law \cite{Bertsekas:DataNetworks}.  Define
\begin{align*}
  \bar{a}_n^r \defeq \lim_{T \rightarrow \infty} \frac{1}{T} \sum_{t = t_0}^{t_0+T-1}\expect{a_n^r(t)},~
  \bar{a}_n \defeq \lim_{T \rightarrow \infty} \frac{1}{T} \sum_{t = t_0}^{t_0+T-1}\expect{ a_n(t)}.
\end{align*}

\begin{lemma}
  \label{lem:delay}
  Suppose $B > 2(\delta^\maxvar+D)$.  In the steady state, the average delay at real queue $n \in \set{N}$ under the floating-queue algorithm is bounded by
\begin{equation*}
  \text{Per-hop delay} \leq \frac{B}{\bar{a}_n^r - \delta^\maxvar c^\ast e^{-[B/2 - \delta^\maxvar -D]/K}}.
\end{equation*}
\end{lemma}
\begin{IEEEproof}
Since the buffer size of queue $n \in \set{N}$ is $B$, Little's law implies:
\begin{align*}
  & \text{Per-hop Delay} = B\biggl/\prts{\lim_{T \rightarrow \infty} \frac{1}{T} \sum_{t = t_0}^{t_0+T-1} \expect{a_n^{r'}(t)} } \\
  & \quad = B\biggl/\prts{\lim_{T \rightarrow \infty} \frac{1}{T} \sum_{t = t_0}^{t_0+T-1} \expect{ a_n^r(t) - d_n(t) }} \\
  & \quad \leq B\biggl/\prts{\bar{a}_n^r - \delta^\maxvar c^\ast e^{-[B/2 -\delta^\maxvar-D]/K}}.
\end{align*}
\end{IEEEproof}

The implication of Lemma \ref{lem:delay} is that, when $B$ is large enough such that $\delta^\maxvar c^\ast e^{-[B/2 - \delta^\maxvar-D]/K}$ and the number of drops at other queues are negligible, $\bar{a}_n^r$ is approximately $\bar{a}_n$, and the average delay is $O(B)$.

\subsubsection{Objective Cost}
The average objective cost is considered in two cases.  Let $f^\fqvar(t)$ denote the cost under the floating-queue algorithm at time $t$, and $f^\fqvar \defeq \lim_{T \rightarrow \infty} \frac{1}{T} \sum_{t = t_0}^{t_0+T-1} \expect{f^\fqvar(t)}$ denote the expected time-average cost under the floating-queue algorithm.

\textbf{Drop-Independent Cost}:\\
In this case, packet drops do not affect the objective cost.  Such cost can be the energy expenditure that is spent to transmit both real and fake packets.  Due to this independence, the average cost follows immediately from the result of drift-plus-penalty policy \eqref{eq:fdpp}.

\begin{theorem}
  \label{thm:independent}
  Suppose each real queue has buffer size $B > 2(\delta^\maxvar+D)$.  When $V > 0$ and packet drops do not incur any penalty cost, the floating-queue algorithm achieves:
\begin{align*}
  f^\fqvar &= f^\dppvar \leq f^\optvar + O(1/V) \\
  \text{Per-hop delay} &\leq O(B/(1 - e^{-B})) = O(B) \\
  \text{Average drops} &\leq O(e^{-B}).
\end{align*}
\end{theorem}
It can be shown that the transient time of the drift-plus-penalty algorithm is $O(V)$, so parameter $V$ cannot be set to infinity.

\textbf{Drop-Dependent Cost}:\\
In this case, packet drops affect the objective cost.  Such cost can be the amount of admitted packets.  Let $\kappa < \infty$ be a maximum penalty cost per one unit of packet drop.  Then we have the following result.

\begin{theorem}
  \label{thm:dependent}
  Suppose $B > 2(\delta^\maxvar+D)$.  When $V > 0$ and $\kappa$ is a maximum penalty cost per one unit of packet drop, the floating-queue algorithm achieves:
\begin{align*}
  f^\fqvar &\leq f^\optvar + O(1/V) + O(e^{-B}) \\
  \text{Per-hop delay} &\leq O(B/(1 - e^{-B})) = O(B) \\
  \text{Average drops} &\leq O(e^{-B})
\end{align*}
\end{theorem}
\begin{IEEEproof}
Recall that $f(t)$ is a cost incurred at time $t$ under the drift-plus-penalty policy.  At each time $t$, we have
\begin{equation*}
  f^\fqvar(t) \leq f(t) + \kappa \sum_{n=1}^N d_n(t).
\end{equation*}
Summing from $t_0$ to $t_0+T-1$, dividing by $T$, and taking an expectation gives
\begin{multline*}
  \frac{1}{T} \sum_{t = t_0}^{t_0+T-1} \expect{f^\fqvar(t)} \leq \frac{1}{T} \sum_{t = t_0}^{t_0+T-1} \expect{ f(t) } \\+ \frac{\kappa}{T} \sum_{t = t_0}^{t_0+T-1} \sum_{n=1}^N \expect{ d_n(t) }.
\end{multline*}
Taking a limit as $T$ approaches infinity gives
\begin{align*}
  & \lim_{T \rightarrow \infty} \frac{1}{T} \sum_{t = t_0}^{t_0+T-1} \expect{f^\fqvar(t)} \\
  & \quad \leq f^\dppvar + \kappa N \delta^\maxvar c^\ast e^{-[B/2 - \delta^\maxvar -D]/K} \\
  & \quad \leq f^\optvar + O(1/V) + O(e^{-B})
\end{align*}
\end{IEEEproof}

\section{Simulation}
\label{sec:simulation}
\begin{figure}
  \centering
  \includegraphics[scale=1.0]{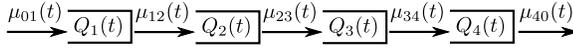}
  \vspace{-0em}
  \caption{Line network}
  \label{fig:simulation}
\end{figure}

A line network with 4 queues, shown in Figure \ref{fig:simulation}, is simulated in two scenarios.  The common network configuration is as follows.  In each time slot, an exogenous packet arrives with probability $0.92$.  Transmission $\mu_{ij}(t)$ is orthogonal and depends on channel state that is ``good'' with probability $0.9$ and ``bad'' with probability $0.1$ for $(i,j) \in \{(1,2), (2,3), (3,4), (4,0)\}$.

\subsection{Power Minimization}
In this scenario, all exogenous arrivals are admitted.  When channel state is ``good'', one packet is transmitted using $1$ unit of power; otherwise $2$ units of power are used.  The goal is to stabilize this network while minimizing the power usage.  Note that the optimal average minimum power is $1 \times 0.9 +2 \times 0.02 = 0.94 $ per hop, and the average total power is $3.76$.

\begin{figure}
  \centering
  \includegraphics[scale=0.43]{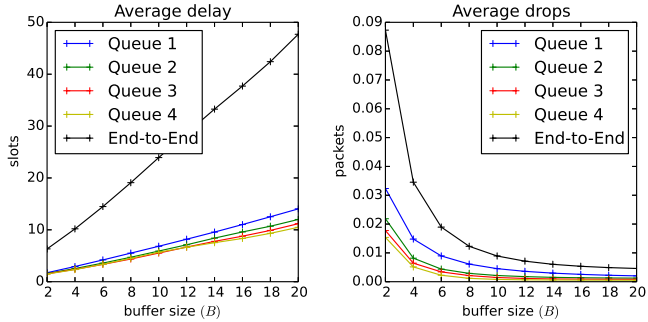}
  \vspace{-1.5em}
  \caption{Results of power minimization problem with $V=200$}
  \label{fig:psim}
\end{figure}

Simulation results of this scenario are shown Figure \ref{fig:psim}.  The time average power expenditure is $3.761$ for all buffer sizes $B$.  In Figure \ref{fig:psim}, the average delay increases linearly with the buffer size, and the average drops decrease exponentially with the buffer size.  This result confirms the bounds in Theorem \ref{thm:independent}.

\subsection{Throughput Maximization}
\begin{figure}
  \centering
  \includegraphics[scale=0.43]{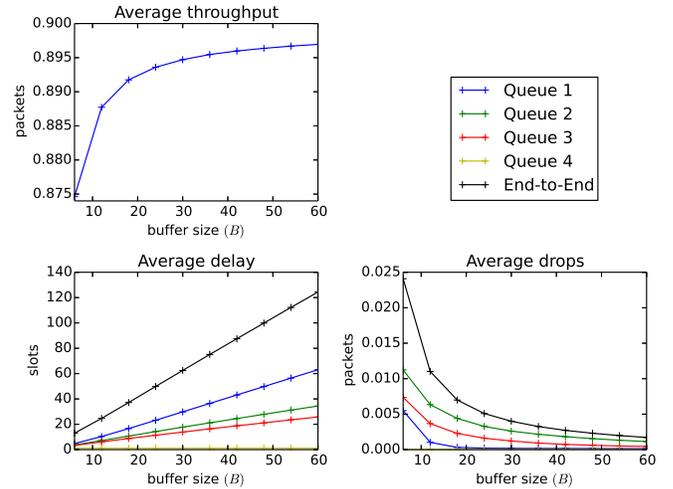}
  \vspace{-1.5em}
  \caption{Results of throughput maximization problem with $V=200$}
  \label{fig:tsim}
\end{figure}

In this scenario, a network decides to admit random exogenous arrival in each time slot.  The goal is to maximize the time-average end-to-end throughput, which are real packets.  Packet drops reduce the value of this objective function.  Transmission $\mu_{ij}(t) = 1$ is possible if its channel state is ``good''; otherwise the transmission is not allowed.    Note that the maximum admission rate is $0.9$, because of the limitation of the average transmission rate.  Figure \ref{fig:tsim} shows the simulation results of this scenario, which comply with the bounds in Theorem \ref{thm:dependent}.

\section{Conclusion}
\label{sec:conclusion}
We propose the general floating-queue algorithm that allows the stochastic network optimization framework to operate with finite buffers.  When the buffer size at each queue is $B$, we prove the proposed algorithm achieves within $O(e^{-B})$ of optimal utility, while the average per-hop delay is $O(B)$.  The finiteness incurs $O(e^{-B})$ drops, decreasing exponentially.  We confirm the theoretical results with simulations.

\bibliographystyle{IEEEtran}
\bibliography{Reference}

\appendix

To prove case 2 of Lemma \ref{lem:TL}, we first show that, for any sample path $\prtc{ \spath_n(t) }_{t=\tau}^{\tau'}$ whose $Q_n^f(t)$ is non-decreasing, there is a simple one that preserves the queue dynamics and the amount of admitted real arrivals under the floating-queue algorithm and the lower-bound policy.  This simple sample path is then used to simplify the later proofs.


\begin{lemma}
  \label{lem:spathconv1}
Given sample path $\prtc{ \spath_n(t) }_{t=\tau}^{\tau'}$ satisfying all dynamics in Section \ref{sec:floatingqueue} where $\spath_n(t) =$
\begin{equation*}
  \prtr{ a_n^r(t), a_n^f(t), a_n^{r'}(t), a_n^{f'}(t), b_n^r(t), b_n^f(t), Q_n^r(t), Q_n^f(t) }
\end{equation*}
with $Q_n^f(t) \leq Q_n^f(t+1)$ for $t \in \prtc{\tau, \dotsc, \tau'}$,
a new sample path $\prtc{ \tilde{\spath}_n(t) }_{t=\tau}^{\tau'}$ where $\tilde{\spath}_n(t) =$
\begin{equation*}
  \prtr{ a_n^r(t), a_n^f(t), a_n^{r'}(t), a_n^{f'}(t), b_n^r(t), \tilde{b}_n^f(t), Q_n^r(t), Q_n^f(t) }
\end{equation*}
with $\tilde{b}_n^f(t) = \min\prts{Q_n^f(t), b_n^f(t)}$ preserves all dynamics in Section \ref{sec:floatingqueue}.
\end{lemma}
\begin{IEEEproof}
From the dynamic of fake queue \eqref{eq:qnf}, $\tilde{b}_n^f(t)$ is the actual number of fake packets that is served at time $t$ under the floating-queue algorithm, i.e., $Q_n^f(t+1) = \max\prts{ Q_n^f(t) - b_n^f(t), 0 } + a_n^{f'}(t) = Q_n^f(t) - \tilde{b}_n^f(t) + a_n^{f'}(t)$.  Therefore, the dynamic of fake queue \eqref{eq:qnf} in the new sample path is still valid.  The dynamic of real queue \eqref{eq:qnr} is valid, because the real part is not modified.


\end{IEEEproof}

The implication of the Lemma \ref{lem:spathconv1} is that, instead of an original sample path, we can consider its alternate version, in which fake services do not excess the amount of fake backlogs.  Both sample paths shares the same admitted real arrivals $a_n^{r'}(t)$ under the floating-queue algorithm.  Further, both sample paths shares the same the admitted real arrivals $\hat{a}_n^r(t)$ under the lower-bound policy, because $Q_n(t)$ and $a_n^r(t)$ of both sample paths are identical.  Notice that we need to preserve all dynamics in Section \ref{sec:floatingqueue}, in order to use those dynamics and properties in Section \ref{sec:floatingqueue} and Section \ref{sec:lowerboundpolicy}, when an alternate sample path is considered.

It is possible to simplify this alternate sample path.

\begin{lemma}
  \label{lem:spathconv2}
Suppose $B \geq 2\delta^\maxvar$. Then given sample path $\prtc{ \spath_n(t) }_{t=\tau}^{\tau'}$ satisfying dynamics in Section \ref{sec:floatingqueue} where $\spath_n(t) = $
\begin{equation*}
  \prtr{ a_n^r(t), a_n^f(t), a_n^{r'}(t), a_n^{f'}(t), b_n^r(t), b_n^f(t), Q_n^r(t), Q_n^f(t) }
\end{equation*}
with $Q_n^f(t) \leq Q_n^f(t+1)$ and $b_n^f(t) \leq Q_n^f(t)$ for $t \in \prtc{\tau, \dotsc, \tau'}$, 
a new sample path $\prtc{ \tilde{\spath}_n(t) }_{t=\tau}^{\tau'}$ where $\tilde{\spath}_n(t) =$
\begin{equation*}
  \prtr{ a_n^r(t), \tilde{a}_n^f(t), a_n^{r'}(t), \tilde{a}_n^{f'}(t), b_n^r(t), \tilde{b}_n^f(t), Q_n^r(t), Q_n^f(t) }
\end{equation*}
with \\
i) $\tilde{b}_n^f(t) = 0$,\\
ii) $\tilde{a}_n^f(t) = a_n^f(t) - b_n^f(t)$, and \\
iii) $\tilde{a}_n^{f'}(t) = \tilde{a}_n^f(t) + (a_n^r(t) - a_n^{r'}(t))$\\
preserves all dynamics in Section \ref{sec:floatingqueue}.
\end{lemma}
\begin{IEEEproof}
For the given sample path $\prtc{ \spath_n(t) }_{t=\tau}^{\tau'}$, fake backlogs are non-decreasing and $Q_n^f(t) \leq Q_n^f(t+1)$ for $t \in \prtc{ \tau, \dotsc, \tau' }$.  Since $b_n^f(t) \leq Q_n^f(t)$, it follows from \eqref{eq:qnf} that $Q_n^f(t+1) = Q_n^f(t) - b_n^f(t) + a_n^{f'}(t)$.  From the non-decreasing, we have
\begin{equation}
  \label{eq:sc21}
  b_n^f(t) \leq a_n^{f'}(t).
\end{equation} 
From \eqref{eq:bnr} and \eqref{eq:bnf}, if $b_n^f(t) > 0$, then $Q_n^r(t) < b_n(t) \leq \delta^\maxvar$ and $d_n(t) = 0$, since $B \geq 2\delta^\maxvar$.  From \eqref{eq:anf}, it follows that $a_n^{f'}(t) = a_n^f(t)$.  Using \eqref{eq:sc21}, we have $b_n^f(t) \leq a_n^f(t)$ and $a_n^f(t) - b_n^f(t) \geq 0$.  Therefore, the new sample path with $\tilde{a}_n^f(t) = a_n^f(t) - b_n^f(t)$ and $\tilde{b}_n^f(t) = 0$ preserves the fake queue dynamic \eqref{eq:qnf}, i.e.,
\begin{align*}
  Q_n^f(t+1) & = \max\prts{Q_n^f(t) - b_n^f(t), 0} + a_n^f(t) + d_n(f)\\
  & = Q_n^f(t) - b_n^f(t) + a_n^f(t) + (a_n^r(t) - a_n^{r'}(t))\\
  & = Q_n^f(t) - \tilde{b}_n^f(t) + \tilde{a}_n^f(t) + (a_n^r(t) - a_n^{r'}(t)) \\
  & = \max\prts{Q_n^f(t) - \tilde{b}_n^f(t), 0} + \tilde{a}_n^{f'}(t),
\end{align*}
where $\tilde{a}_n^{f'}(t) = \tilde{a}_n^f(t) + (a_n^r(t) - a_n^{r'}(t))$.  The dynamic of real queue \eqref{eq:qnr} is valid, because the real part is not modified.
\end{IEEEproof}

The implication of Lemma \ref{lem:spathconv1} and Lemma \ref{lem:spathconv2} is as follows.  Given an original sample path $\prtc{ \spath_n(t) }_{t=\tau}^{\tau'}$ with non-decreasing fake backlogs, we can consider an alternate simple sample path $\prtc{ \tilde{\spath}_n(t) }_{t=\tau}^{\tau'}$, in which all dynamics and properties of the floating-queue algorithm and the lower-bound policy still holds.  Technically, the fake-backlog dynamic \eqref{eq:qnf} becomes simpler, i.e., $Q_n^f(t+1) = Q_n^f(t) + \tilde{a}_n^f(t)$ for all $t \in \prtc{\tau, \dotsc, \tau'}$ because $\tilde{b}_n^f(t) = 0$.
Also, the amount of admitted real arrivals under both float-queue algorithm and lower-bound policy are not changed.
These simplify the later proofs.  Before the main lemma, several properties are proved.


\begin{lemma}
  \label{lem:tkj}
  Suppose $B \geq 2\delta^\maxvar$.  Given any sample path $\prtc{ \spath_n(t) }_{t=t_0}^{t_0+T}$ satisfying all dynamics in Section \ref{sec:floatingqueue} for any $t_0$ and positive integer $T$, the following holds for $k \in \prtc{1, \dotsc, K(T)}$, $j \in \prtc{1, \dotsc, J(k)}$, $t_{1,1}^- - 1 \geq t_0$, and $t_{K(T),J(K(T))}^+ +1 \leq t_0+T$:
\begin{itemize}
\item $Q_n^r(t_{k,j}^-) = a_n^r(t_{k,j}^- -1)$,
\item $\hat{a}_n^r(t_{k,j}^- - 1) = 0$,
\item $Q_n^f(t_{k,j}^- - 1) < \dvar_n^{V\ast} - B/2 + \delta^\maxvar$,
\item $Q_n^f(t_{k,j}^+ + 1) < \dvar_n^{V\ast} - B/2 + \delta^\maxvar$.
\end{itemize}
\end{lemma}
\begin{IEEEproof}
For the first property, all real backlogs $Q_n^r(t_{k,j}^- -1)$ must be served at time $t_{k,j}^- -1$, in order to have $b_n^f(t_{k,j}^- -1) > 0$ and $Q_n^f(t_{k,j}^- -1) > Q_n^f(t_{k,j}^-)$, so $Q_n^r(t_{k,j}^- -1) - b_n^r(t_{k,j}^- -1) = 0$.  From Lemma \ref{lem:dropincrease}, $Q_n^f(t_{k,j}^-) < Q_n^f(t_{k,j}^- -1)$ implies $a_n^r(t) = a_n^{r'}(t)$.  Thus, it follows from \eqref{eq:qnr} that $Q_n^r(t_{k,j}^-) = 0 + a_n^r(t_{k,j}^- -1)$, which proves the first property.

For the second property, Let $\tilde{b}_n(t_{k,j}^- -1 ) = \min\prts{ Q_n(t_{k,j}^- -1), b_n(t_{k,j}^- -1)}$.  We shows that $Q_n(t_{k,j}^- -1) \notin \set{B}_n$ by considering the dynamic of the standard queue \eqref{eq:queue_simple}:
\begin{align*}
  Q_n( t_{k,j}^- - 1)
  & = Q_n(t_{k,j}^-) + \tilde{b}_n(t_{k,j}^--1) - a_n(t_{k,j}^- - 1) \\
  & \hspace{-4em}= Q_n^r(t_{k,j}^-) + Q_n^f( t_{k,j}^-) + \tilde{b}_n(t_{k,j}^--1) - a_n(t_{k,j}^--1) \\
  & \hspace{-4em}= a_n^r(t_{k,j}^- - 1) + Q_n^f( t_{k,j}^-) + \tilde{b}_n(t_{k,j}^--1) - a_n(t_{k,j}^--1) \\
  & \hspace{-4em}< \dvar_n^{V\ast} - B/2 + \delta^\maxvar \notin \set{B}_n.
\end{align*}
The last equality uses the first property.  The last inequality uses the facts that $a_n^r(t_{k,j}^- -1) \leq a_n(t_{k,j}^- -1)$ and $Q_n^f( t_{k,j}^- ) < \dvar_n^{V\ast} - B/2$.

The third property can be proved from \eqref{eq:qnf} as follows:
\begin{align*}
  Q_n^f(t_{k,j}^-) & = \max\prts{ Q_n^f(t_{k,j}^- -1) - b_n^f(t_{k,j}^- -1), 0} + a_n^{f'}(t_{k,j}^- -1) \\
  & \geq Q_n^f(t_{k,j}^- -1) - b_n^{f}(t_{k,j}^- -1),
\end{align*}
and
\begin{equation*}
  Q_n^f(t_{k,j}^- -1) \leq Q_n^f(t_{k,j}^-) + b_n^f(t_{k,j}^- -1) \leq \dvar_n^{V\ast} - B/2 + \delta^\maxvar.
\end{equation*}
The last inequality uses the fact that $Q_n^f(t_{k,j}^-) < \dvar_n^{V\ast} - B/2$.

Similarly, the last property can be proved from \eqref{eq:qnf} as follows:
\begin{align*}
  Q_n^f(t_{k,j}^+ + 1) & = \max\prts{ Q_n^f(t_{k,j}^+) - b_n^f(t_{k,j}^+), 0} + a_n^{f'}(t_{k,j}^+) \\
  & \leq Q_n^f(t_{k,j}^+) + a_n^{f'}(t_{k,j}^+) \\
  & < \dvar_n^{V\ast} - B/2 + \delta^\maxvar.
\end{align*}
The last inequality uses the fact that $Q_n^f(t_{k,j}^+) < \dvar_n^{V\ast} - B/2$.
\end{IEEEproof}

We now quantify packet drops under the floating queue algorithm via a modified sample path.  The sample paths between $t_{k,j}^-$ and $t_{k,j}^+$ for some $k$ and $j$ can be either ``non-decreasing (at the end)'' type or ``decreasing (at the end)'' type.  This depends on the non-decreasing or decreasing of fake backlogs at the end of the sample path.  For example, in Figure \ref{fig:TL}, a sample path between $t_{2,2}^-$ and $t_{2,2}^+$ is the non-decreasing type, and a sample path between $t_{2,1}^-$ and $t_{2,1}^+$ is the decreasing type.  We prove each type separately.

\begin{lemma}
  \label{lem:case2_nd}
  When the buffer size of real queue $n$ is $B \geq 2\delta^\maxvar$, given a non-decreasing-type sample path $\prtc{\spath_n(t)}_{t=t_{k,j}^- -1}^{t_{k,j}^+}$ with $b_n^f(t) = 0$ for all $t \in \prtc{t_{k,j}^-, \dotsc, t_{k,j}^+}$,
the following holds
\begin{equation}
  \sum_{t = t_{k,j}^- -1}^{t_{k,j}^+} a_n^{r'}(t) \geq \sum_{t = t_{k,j}^- -1}^{t_{k,j}^+} \hat{a}_n^r(t).
\end{equation}
\end{lemma}
\begin{IEEEproof}
We first consider drops under the floating-queue algorithm.  From \eqref{eq:qnf} with $b_n^f(t) = 0$, it follows that $Q_n^f(t+1) = Q_n^f(t) + a_n^{f'}(t)$ for all $t \in \prtc{ t_{k,j}^-, \dotsc, t_{k,j}^+}$.  Summing from $t = t_{k,j}^-$ to $t_{k,j}^+$ yields
\begin{align*}
  Q_n^f( t_{k,j}^+ +1 ) 
  & = Q_n^f( t_{k,j}^-) + \sum_{t = t_{k,j}^-}^{t_{k,j}^+} a_n^{f'}(t) \\
  & = Q_n^f( t_{k,j}^-) + \sum_{t = t_{k,j}^-}^{t_{k,j}^+} \prts{ a_n^{f}(t) + d_n(t) },
\end{align*}
where the last equality uses \eqref{eq:anf}.  Rearranging terms gives
\begin{equation*}
  \sum_{t = t_{k,j}^-}^{t_{k,j}^+} d_n(t) = Q_n^f( t_{k,j}^+ +1 ) - Q_n^f( t_{k,j}^- ) - \sum_{t = t_{k,j}^-}^{t_{k,j}^+ -1} a_n^f(t).
\end{equation*}
Since $Q_n^f(t_{k,j}^-)$ is a local minimum and $Q_n^f( t_{k,j}^- ) < Q_n^f( t_{k,j}^- -1 )$, Lemma \ref{lem:dropincrease} implies that $d_n( t_{k,j}^- -1 ) = 0$.  From the above equation, we have that
\begin{equation}
  \label{eq:case2_nd_dfp}
  \sum_{t = t_{k,j}^- -1}^{t_{k,j}^+} d_n(t) = Q_n^f( t_{k,j}^+ +1 ) - Q_n^f( t_{k,j}^- ) - \sum_{t = t_{k,j}^-}^{t_{k,j}^+} a_n^f(t).
\end{equation}

Now, drops under the lower-bound policy is considered.  We suppose that $\sum_{t = t_{k,j}^- -1}^{t_{k,j}^+} d_n(t) > 0$; otherwise, there is no drop under the floating-queue algorithm, and
\begin{equation*}
  \sum_{t = t_{k,j}^- -1}^{t_{k,j}^+} a_n^{r'}(t) = \sum_{t = t_{k,j}^- -1}^{t_{k,j}^+ } a_n^{r}(t) \geq \sum_{t = t_{k,j}^- -1}^{t_{k,j}^+ } \hat{a}_n^{r}(t).
\end{equation*}

From Lemma \ref{lem:tkj}, we know that $\hat{a}_n^r(t_{k,j}^- -1) = 0$.  Therefore, we suppose that $Q_n(t) \geq \dvar_n^{V\ast} - B/2 + \delta^\maxvar$ for some $t \in \prtc{ t_{k,j}^-, \dotsc, t_{k,j}^+}$; otherwise, the lower-bound policy drops all real arrivals, and 
\begin{equation*}
  \sum_{t = t_{k,j}^- -1}^{t_{k,j}^+} \hat{a}_n^{r}(t) = 0 \leq \sum_{t = t_{k,j}^- -1}^{t_{k,j}^+} a_n^{r'}(t).
\end{equation*}

When the above assumptions are imposed, we let
\begin{equation*}
  t_{k,j}^* = \arginf_{t \in \prtc{t_{k,j}^-, \dotsc, t_{k,j}^+}} \prtc{ Q_n(t) \geq \dvar_n^{V\ast} - B/2 + \delta^\maxvar },
\end{equation*}
be the first time that $Q_n(t)$ is at least $\dvar_n^{V\ast} - B/2 + \delta^\maxvar$.  This is illustrated in Figure \ref{fig:TL}.  Notice that $t_{k,j}^* \in \prtc{ t_{k,j}^-, \dotsc, t_{k,j}^+}$ by the later assumption.

From \eqref{eq:qnr} and \eqref{eq:qnf} with $b_n^f(t) = 0$, it holds for $t \in \prtc{t_{k,j}^-, \dotsc, t_{k,j}^* -1}$ that
\begin{align*}
  Q_n(t+1) 
  & = Q_n^r(t+1) + Q_n^f(t+1) \\
  & = Q_n^r(t) - b_n^r(t) + a_n^{r'}(t) + Q_n^f(t) + a_n^{f'}(t) \\
  & = Q_n(t) - b_n^r(t) + a_n^{r'}(t) + a_n^{f'}(t).
\end{align*}

Summing the above equation from $t = t_{k,j}^-$ to $t_{k,j}^* -1$ yields
\begin{align*}
  Q_n( t_{k,j}^* )
  & = Q_n(t_{k,j}^-) + \sum_{t = t_{k,j}^-}^{t_{k,j}^* -1} \prts{ -b_n^r(t) + a_n^{r'}(t) + a_n^{f'}(t) } \\
  & = Q_n(t_{k,j}^-) + \sum_{t = t_{k,j}^-}^{t_{k,j}^* -1} \prts{ -b_n^r(t) + a_n^{r}(t) + a_n^{f}(t) },
\end{align*}
and
\begin{align*}
  \sum_{t = t_{k,j}^-}^{t_{k,j}^* -1} a_n^r(t)
  & = Q_n(t_{k,j}^*) - Q_n(t_{k,j}^-) + \sum_{t = t_{k,j}^-}^{t_{k,j}^* -1} \prts{ b_n^r(t) - a_n^{f}(t) } \\
  & \hspace{-4em}\geq Q_n(t_{k,j}^*) - Q_n(t_{k,j}^-) - \sum_{t = t_{k,j}^-}^{t_{k,j}^* -1} a_n^{f}(t)  \\
  & \hspace{-4em}\geq \dvar_n^{V\ast} - B/2 + \delta^\maxvar - Q_n(t_{k,j}^-) - \sum_{t = t_{k,j}^-}^{t_{k,j}^* -1} a_n^{f}(t) \\
  & \hspace{-4em}> Q_n^f(t_{k,j}^+ +1) - \prts{ a_n^r(t_{k,j}^- -1) + Q_n^f(t_{k,j}^-) } - \sum_{t = t_{k,j}^-}^{t_{k,j}^* -1} a_n^{f}(t).
\end{align*}
The last steps uses Lemma \ref{lem:tkj} and the facts that $Q_n(t_{k,j}^-) = Q_n^r(t_{k,j}^-) + Q_n^f(t_{k,j}^-)$.

From the definition of $t_{k,j}^*$, we know that $Q_n(t) < \dvar_n^{V\ast} - B/2 + \delta^\maxvar$ and all real arrivals under the lower-bound policy are dropped, $a_n^r(t) = \hat{d}_n(t)$, for all $t \in \prtc{ t_{k,j}^-, \dotsc, t_{k,j}^*-1 }$.  The third property in Lemma \ref{lem:tkj} also implies that $a_n^r(t_{k,j}^- -1) = \hat{d}_n(t_{k,j}^- -1)$.  Therefore, the above inequality yields
\begin{align*}
  \sum_{t = t_{k,j}^- -1}^{t_{k,j}^* -1} \hat{d}_n(t) &= \sum_{t = t_{k,j}^- -1}^{t_{k,j}^* -1} a_n^r(t) \\
  & > Q_n^f(t_{k,j}^+ +1) - Q_n^f(t_{k,j}^-) - \sum_{t = t_{k,j}^-}^{t_{k,j}^* -1} a_n^{f}(t) \\
  & \geq Q_n^f(t_{k,j}^+ +1) - Q_n^f(t_{k,j}^-) - \sum_{t = t_{k,j}^-}^{t_{k,j}^+ } a_n^{f}(t) \\
  & = \sum_{t = t_{k,j}^- -1}^{t_{k,j}^+} d_n(t),
\end{align*}
where the last step uses \eqref{eq:case2_nd_dfp}.  Then we have
\begin{equation*}
  \sum_{t = t_{k,j}^- -1}^{t_{k,j}^+} \hat{d}_n(t) > \sum_{t = t_{k,j}^- -1}^{t_{k,j}^+} d_n(t).
\end{equation*}
Applying \eqref{eq:dn} and \eqref{eq:hdn} proves the lemma.  Note that the equality holds in some situations when the both algorithms do not drop any real packets.
\end{IEEEproof}

\begin{lemma}
  \label{lem:case2_d}
  When the buffer size of real queue $n$ is $B \geq 2\delta^\maxvar$, given a decreasing-type sample path $\prtc{\spath_n(t)}_{t=t_{k,j}^- -1}^{t_{k,j}^+}$ with $b_n^f(t) = 0$ for all $t \in \prtc{t_{k,j}^-, \dotsc, t_{k,j}^+ -1}$ and $Q_n^f(t_{k,j}^+) > Q_n^f(t_{k,j}^+ +1)$,
the following holds
\begin{equation*}
  \sum_{t = t_{k,j}^- -1}^{t_{k,j}^+} a_n^{r'}(t) \geq \sum_{t = t_{k,j}^- -1}^{t_{k,j}^+} \hat{a}_n^r(t).
\end{equation*}
\end{lemma}
\begin{IEEEproof}
We construct a new sample path $\prtc{\tilde{\spath}_n(t)}_{t=t_{k,j}^-}^{t_{k,j}^+ -1}$ where $\tilde{\spath}_n(t) = \spath_n(t)$ for all $t \in \prtc{t_{k,j}^-, \dotsc, t_{k,j}^+ -1}$.  Applying Lemma \ref{lem:case2_nd} on the new sample path implies that
\begin{equation*}
  \sum_{t = t_{k,j}^- -1}^{t_{k,j}^+ -1} a_n^{r'}(t) \geq \sum_{t = t_{k,j}^- -1}^{t_{k,j}^+ -1} \hat{a}_n^r(t).
\end{equation*}

Since $Q_n^f(t_{k,j}^+ +1) < Q_n^f(t_{k,j}^+)$, Lemma \ref{lem:dropincrease} implies that $a_n^{r'}(t_{k,j}^+) = a_n^r(t_{k,j}^+)$.  This means that $a_n^{r'}(t_{k,j}^+) \geq \hat{a}_n^r(t_{k,j}^+)$.  Using this and the above relation proves the lemma.
\end{IEEEproof}

Finally, the second case of Lemma \ref{lem:TL} can be proved.
\begin{lemma}
  \label{lem:case2}
  When the buffer size of real queue $n \in \set{N}$ is $B \geq 2\delta^\maxvar$, given any sample path $\prtc{ \spath_n(t) }_{t=t_{k,j}^- -1}^{t=t_{k,j}^+}$ satisfying all dynamics in Section \ref{sec:floatingqueue}, it holds that
\begin{equation*}
  \sum_{t = t_{k,j}^- -1}^{t_{k,j}^+} a_n^{r'}(t) \geq \sum_{t = t_{k,j}^- -1}^{t_{k,j}^+} \hat{a}_n^r(t).
\end{equation*}
\end{lemma}
\begin{IEEEproof}
We construct a new sample path and apply a previous lemma, depending on the type of the original path.  Let $\prtc{ \tilde{\spath}_n(t) }_{t=t_{k,j}^- -1}^{t_{k,j}^+}$ be a new sample path with $\tilde{\spath}_n(t_{k,j}^- -1) = \spath_n(t_{k,j}^- -1)$.

When $Q_n^f(t_{k,j}^+ + 1) \geq Q_n^f(t_{k,j}^+)$ (non-decreasing type), let $\prtc{ \tilde{\spath}_n(t) }_{t=t_{k,j}^-}^{t_{k,j}^+}$ be a result of applying Lemma \ref{lem:spathconv1} and Lemma \ref{lem:spathconv2} on $\prtc{ \spath_n(t) }_{t=t_{k,j}^-}^{t_{k,j}^+}$.  Then, applying Lemma \ref{lem:case2_nd} on the new sample path proves the lemma.

When $Q_n^f(t_{k,j}^+ + 1) < Q_n^f(t_{k,j}^+)$ (decreasing type), let $\prtc{ \tilde{\spath}_n(t) }_{t=t_{k,j}^-}^{t_{k,j}^+ -1}$ be a result of applying Lemma \ref{lem:spathconv1} and Lemma \ref{lem:spathconv2} on $\prtc{ \spath_n(t) }_{t=t_{k,j}^-}^{t_{k,j}^+ -1}$, and let $\tilde{\spath}_n(t_{k,j}^+) = \spath_n(t_{k,j}^+)$.  Then, applying Lemma \ref{lem:case2_d} on the new sample path proves the lemma.
\end{IEEEproof}

\end{document}